\documentclass[MSNbibl,nameyear,dvips]{arxbj}

\aid{0}
\volume{19}
\issue{4}
\pubyear{2013}
\firstpage{1449}
\lastpage{1464}
\doi{10.3150/12-BEJSP15} 

\makeatletter

\define@key{bid}{pmid}{}
\define@key{bid}{pmcid}{}
\define@key{bid}{mid}{}

\newcommand{\rrvert}{\vert}
\newcommand{\llvert}{\vert}

\newcommand{\cal}{\mathcal}

\newcommand{\implies}{\Longrightarrow}

\newcommand{\argmax}{\mathop{\operatorname{argmax}}}
\newcommand{\argmin}{\mathop{\operatorname{argmin}}}
\newcommand{\minimize}{\mathop{\operatorname{minimize}}}
\newcommand{\maximize}{\mathop{\operatorname{maximize}}}

\newcommand{\real}{\mathbb{R}}
\newcommand{\Pp}{{\mathbb P}}
\newcommand{\Ee}{{\mathbb E}}
\newcommand{\lips}{{\cal L}}
\newcommand{\tube}{\operatorname{Tube}}
\newcommand{\XK}{K_X}

\makeatother

\begin{document}
\begin{frontmatter}

\title{The geometry of least squares in the 21st century}
\runtitle{Method of least squares}

\begin{aug}
\author{\fnms{Jonathan} \snm{Taylor}\corref{}\ead[label=e1]{jonathan.taylor@stanford.edu}}
\runauthor{J. Taylor} 
\address{Department of Statistics, Stanford University, Sequoia
Hall, 390 Serra Mall, Stanford, CA 94305, USA. \printead{e1}}
\end{aug}


%
\begin{abstract}
It has been over 200 years since Gauss's and Legendre's famous priority
dispute on who discovered the method of least squares. Nevertheless, we
argue that the normal equations are still relevant in many facets of
modern statistics, particularly in the domain of high-dimensional
inference. Even today, we are still learning new things about the law
of large numbers, first described in Bernoulli's \textit{Ars Conjectandi}
300 years ago, as it applies to high dimensional inference.

The other insight the normal equations provide is the asymptotic
Gaussianity of the least squares estimators. The general form of the
Gaussian distribution, Gaussian processes, are another tool used in
modern high-dimensional inference. The Gaussian distribution also
arises via the central limit theorem in describing weak convergence of
the usual least squares estimators. In terms of high-dimensional
inference, we are still missing the right notion of weak convergence.

In this mostly expository work, we try to describe how both the normal
equations and the theory of Gaussian processes, what we refer to as the
``geometry of least squares,'' apply to many questions of current
interest.
\end{abstract}

%

\end{frontmatter}

\section{Basic tools in the geometry of least squares}\label{secintro}

The method of least squares has by now a long and well-trodden history,
which we will not attempt
to address in this work. Toward the end of the 20th century,
\citet{stigler} referred to the method of least squares as the automobile
of (then) modern statistical analysis. Today, 30 years later, as
automobiles have
modernized, including technological and efficiency improvements, so too
have the methods of least
squares changed.

Let us recall the classical least squares problem: given outcome vector
$y \in\real^n$ and design matrix
$X \in\real^{n \times p}$, the least squares problem is typically
posed as
%
\begin{equation}
\label{eqLS} \minimize_{\beta\in\real^p} \frac{1}{2} \|y-X\beta
\|^2_2.
\end{equation}
The related normal equations are
%
\begin{equation}
\label{eqLSKKT} X^T(X\hat{\beta}-y) = 0.
\end{equation}

We could have almost equivalently written the equations above in terms of
the metric projection problem
%
\begin{equation}
\label{eqmetricproj} y \mapsto\argmin_{x \in A} \frac{1}{2} \|y-x
\|^2_2 \stackrel{\Delta} {=} \pi_A(y)
\end{equation}
with classical least squares by setting $A=\operatorname{col}(X)$ the
column space of $X$. This metric projection problem
is the first of the basic tools in what we will refer to here as ``the
geometry of least squares in the 21st century.''

While written as a function above, the map (\ref{eqmetricproj}) is
not always a function, there may be many minimizers for a given $y$.
For a set $A$, define the \textit{critical radius} of $A$ as
%
\begin{equation}
\label{eqradius} r_c(A) = \sup\Bigl\{r: \inf_{x \in A}
\|y-x\|_2 \leq r \implies\mbox{(\ref{eqmetricproj}) has a unique
solution} \Bigr\}.
\end{equation}
Note that for convex $A$, $r_c(A) = +\infty$. The metric projection
problem also
makes sense for other sets and other metrics. For instance, suppose $A
\subset S(\real^n)=S_{\ell_2}(\real^n)$, the unit
$\ell_2$ sphere in $\real^n$ with distance
$d(x,y) = \cos^{-1}(x^Ty)$. One might also consider the spherical
metric projection
%
\begin{equation}
\label{eqmetricprojsphere} S\bigl(\real^n\bigr) \ni y \mapsto
\argmax_{x \in A} x^Ty
\end{equation}
with the critical radius (\ref{eqradius}) being similarly defined. We
will see in Section~\ref{secgaussian} that the above critical radius
plays a part in the supremum of a class of Gaussian processes
[\citet{rfg}], one of the other important class of objects
associated to
Gauss's name. Gaussian processes suffer some of the same deficiencies
identified in \citet{stigler}: they make many assumptions and have their
limitations. Nevertheless, they are a crucial inferential tool in
analyzing the behavior of (\ref{eqmetricproj}). Hence,
we refer to these as the second of the basic tools in the geometry of
least squares in the 21st century.

\section{A canonical high-dimensional regression problem}
\label{secproblem}

In the classical setting $n > p$ the system (\ref{eqLSKKT}) often
has a unique solution, the familiar
\[
\hat{\beta} = \bigl(X^TX\bigr)^{-1}\bigl(X^Ty
\bigr).
\]
In many parametric models, the least squares model is of course too
simple. In the exponential family setting
[\citet{amari}, \citet{efron}], the normal equations are similar, with $(X^TX)$
replaced by the observed Fisher information. We have focused on squared
error-loss
for its simplicity of exposition.

In high-dimensional settings, $n$ is often less than $p$ and there is
of course no unique solution
to (\ref{eqLSKKT}). Many modifications are possible, for instance,
ridge or Tikhonov regularization which adds a strongly convex quadratic
term to (\ref{eqLS}). The addition of such quadratic terms changes
the quadratic part of the loss but does not fundamentally\vadjust{\goodbreak} change much
else until we begin
to make assumptions about whether or not the model is correct, and how
much bias might be incurred
by such regularization.

In modern high-dimensional settings the regularization term, or
penalty, of choice is often a norm,
with the LASSO [\citet{tibshiranilasso}] being the most popular. The
lasso problem is
%
\begin{equation}
\label{eqlasso0} \minimize_{\beta\in\real^p} \frac{1}{2} \|y-X\beta
\|^2_2 + \lambda\|\beta\|_1.
\end{equation}
The duality between norms allows us to write
\[
\|\beta\|_1 = \sup_{u:\|u\|_{\infty} \leq1} u^T\beta=
h_{B_{\infty
}}(\beta)
\]
with $B_{\infty}$ the $\ell_{\infty}$ ball of radius 1 in $\real^p$ and
for any set $K$
%
\begin{equation}
\label{eqsupport} h_K(\beta) = \sup_{u \in K}
u^T\beta
\end{equation}
is the support function of the set $K$ which we assume to be closed and
containing 0. In this notation,
the LASSO problem can be expressed as
%
\begin{equation}
\label{eqlasso} \minimize_{\beta\in\real^p} \frac{1}{2} \|y-X\beta
\|^2_2 + \lambda h_{B_{\infty}}(\beta).
\end{equation}

Our canonical problem is therefore
%
\begin{equation}
\label{eqcanonical}
\minimize_{\beta\in\real^p}
\frac{1}{2} \|y-X\beta\|^2_2+ \lambda
h_{K}(\beta)
\end{equation}
with $K$ being a closed, convex set containing 0. For one of many
possible infinite
dimensional formulations of this canonical problem, see \citet
{tsirelson} whose author is also associated
to one of the most famous tools in the theory of Gaussian process, the
Borell TIS inequality
[\citet{rfg}].

The normal equations of the least squares problem are replaced with the
KKT conditions
[\citet{boyd}] for (\ref{eqcanonical}). For our canonical problem,
the KKT conditions are
%
\begin{equation}
\label{eqKKT} X^T(X\hat{\beta}-y) + \hat{u}= 0, \qquad\hat{u} \in
\lambda\cdot\partial h_K(\hat{\beta}),
\end{equation}
where the $\partial$ denotes the sub differential. In what follows, we denote
a solution to this problem as $\hat{\beta}_{\lambda,\XK}$ to denote
the dependence on
the penalty $\XK$ and the penalty parameter $\lambda$.\looseness=-1

As we can encode linear or cone constraints in the support function,
it is safe to say that a huge number of problems fit into this framework.
Some examples include:
\begin{itemize}
\item LASSO [\citet{tibshiranilasso}];
\item compressed sensing [\citet{candescs}, \citet{donohocs}];
\item group LASSO [\citet{grouplasso}, \citet{overlapgrouplasso}];
\item graphical LASSO [\citet{graphicallasso}, \citet{buhlmann}];
\item matrix completion [\citet{candesrechtexact}, \citet{mazumderhastiesoftsvd}];\vadjust{\goodbreak}
\item sign restricted regression [\citet{meinshausensign}];
\item hierarchically constrained models [\citet
{bienhierarchical}, \citet{bachhierarchical}];
\item generalized LASSO [\citet{tibshiranitaylorpath}];
\item total variation denoising [\citet{nesta}].
\end{itemize}

The relation between (\ref{eqcanonical}) and the metric projection
map is through a particular
dual function, also developed by Legendre. The canonical dual problem is
%
\begin{equation}
\label{eqcanonicaldual} \minimize_{u \in\lambda\cdot\operatorname
{row}(X) \cap K} \frac
{1}{2}
(u-w)^T\bigl(X^TX\bigr)^{\dagger}(u-w)
\end{equation}
with $\operatorname{row}(X)$ the row space of the matrix $X$. This
dual problem can be derived by minimizing
the following Lagrangian with respect to $\beta, \eta$
%
\begin{equation}
\label{eqlagrangian} L(\beta,\eta;u) = \tfrac{1}{2} \|y-X\beta
\|^2_2 + \lambda h_K(\eta) + u^T(
\beta-\eta).
\end{equation}
After a sign change, the problem in (\ref{eqcanonicaldual}) is
fairly easily seen to be equivalent to
%
\begin{equation}
\label{eqcanonicaldual0} \maximize_{u} \Bigl[\inf
_{\beta,\eta} L(\beta,\eta;u) \Bigr]
\end{equation}
with the constraints in (\ref{eqcanonicaldual}) encoding the fact that
\[
\inf_{\beta,\eta} L(\beta,\eta;u) = -\infty
\]
whenever $u \notin\lambda\cdot\operatorname{row}(X) \cap K. $
Any pair $\hat{\beta}, \hat{u}$ is related through
%
\begin{equation}
\label{eqprimaldual} \bigl(X^TX\bigr)\hat{\beta} - w + \hat{u} = 0.
\end{equation}

Choosing an orthonormal basis for $\operatorname{row}(X)$, we see that
the dual problem can be phrased as the metric projection problem of
$(X^TX)^{-\dagger/2}w$ onto
$ \lambda\cdot\operatorname{row}(X) \cap K$. Alternatively, if we
are interested only in
the fitted values $X\hat{\beta}_{\lambda, K}$, the original problem
(\ref{eqcanonical}) can be expressed as the residual $\hat{\mu
}=y-\hat{r}$
from
%
\begin{equation}
\label{eqresid0} \hat{r} = \argmin_{r \in\lambda\cdot\XK} \frac{1}{2} \|y-r
\|^2_2,
\end{equation}
where
%
\begin{equation}
\label{eqXK} \XK= \bigl(X^T\bigr)^{-1}K.
\end{equation}
Or, in another form,
%
\begin{equation}
\label{eqresid} \hat{\mu}(y) = y - \pi_{\lambda\cdot\XK}(y).
\end{equation}
We see that our canonical regression problem is in fact a metric
projection problem.

Having posed the canonical high-dimensional regression problem as a
metric projection problem,
we now try to describe how metric projection is related to some
fundamental issues in the understanding
of this problem both from an algorithmic view and an inferential view.

\subsection{Algorithms to solve the canonical problem}
\label{secalgorithm}

The general problem is phrased in terms of an arbitrary $K$. For many
high-dimensional
problems, this $K$ is chosen to emphasize expected structure in the
data. For example,
it is well known that the LASSO yields sparse solutions, the group
LASSO yields
groups of nonzero coefficients, etc.

This special structure is based on a particular structure encoded in
$K$. Further, for many canonical
choices of $K$ used in high dimensional statistics, such as those cited
in Section~\ref{secproblem},
the following metric projection is simple
\[
\nu\mapsto\argmin_{ u \in\lambda\cdot K} \frac{1}{2} \|\nu-u\|
^2_2 = \nu- \argmin_{\beta} \biggl[
\frac{1}{2} \|\nu-\beta\|^2_2 + \lambda
h_K(\beta) \biggr].
\]

Such optimization problems are referred to as problems in composite form
[Becker, Bobin and Cand{\`e}s (\citeyear{nesta}), \citet{boyd}, \citet{nesterov}]. For such problems,
many modern first order solvers use a version of
generalized gradient descent. The steps in generalized gradient descent
are essentially iterations of this metric projection map. Specifically, to
solve the canonical problem (\ref{eqcanonical}) for a given step-size
$\alpha$ a simple generalized
gradient algorithm reads as
%
\begin{eqnarray}
\hat{\beta}^{(k+1)} &=& \argmin_{\beta}
\frac{1}{2 \alpha} \bigl\|\nu^{(k)}_{\alpha}-\beta\bigr\|^2_2
+ \lambda h_K(\beta)
\nonumber\\
&=& \nu^{(k)}_{\alpha} - \argmin_{\eta\in\lambda\alpha\cdot K} \frac{1}{2}
\bigl\|\nu^{(k)}_{\alpha}-\eta\bigr\|^2_2,
\\
\nu^{(k)}_{\alpha} &=& \beta^{(k)} - \alpha\cdot
X^T\bigl(X\beta^{(k)}-y\bigr).
\nonumber
\end{eqnarray}

The first line in the update above is the usual form of updates for
generalized gradient descent, while the second line
expresses this step as the residual after an application of the metric
projection map.
Accelerated schemes can do much better with slightly different updates
above, see
\citet{nesta}, \citet{nesterov}, \citet{tseng}.
With modern computing techniques, such simple algorithms can scale to huge
problems, see \citet{admm}, \citet{mazumderhastiesoftsvd}.

\section{Inference for the canonical problem}
\label{secinference}

\subsection{Law of large numbers}

Having solved the canonical problem, what can we say about its solution?
As this special issue is devoted to the appearance of one of the first
proofs of the law of large numbers,
we should at least hope to provide such an answer.

In the classical setting, assuming independence, and the usual linear
regression model
%
\begin{equation}
\label{eqmodel} y = \mu+ \varepsilon
\end{equation}
with noise $\varepsilon$ having scale $\sigma$,
the central limit theorem
can often be applied to (\ref{eqLS}) yielding the usual result
%
\begin{equation}
\bigl\|X(\hat{\beta} - \beta_0)\bigr\|_2^2
\stackrel{D} {\approxeq} \sigma^2 \cdot\chi^2_p
\end{equation}
under the null $H_0: \mu=X\beta_0 \in\operatorname{col}(X)$. Of
course, this
forms the basis of much
inference in modern (and not so modern) applied statistics in the fixed
$p$, $n$ growing regime.
In terms of the parameters themselves, this
implies the weaker statement
%
\begin{equation}
\label{eqLLN} \|\hat{\beta} - \beta\|_2 \leq n^{-1/2}\sigma
\cdot\zeta_n
\end{equation}
for some random variable $\zeta_n=O_{\Pp}(1)$.

In the classical setting, assuming $X$ is full rank, the bound (\ref
{eqLLN}) is a simple two line proof followed by
some assertions. If
we write ${\cal L}(\beta) = \frac{1}{2} \|y-X\beta\|^2_2$, then
\begin{eqnarray*}
0 &\geq& {\cal L}(\hat{\beta}) - {\cal L}(
\beta_0)
\\
& = &\nabla{\cal L}({\beta}_0)^T(\hat{\beta}-
\beta_0) + \frac
{1}{2}(\hat{\beta}-\beta_0)^T
\bigl(X^TX\bigr) (\hat{\beta}-\beta_0)
\\
&= &\bigl(X^T\varepsilon\bigr)^T(\hat{\beta}-
\beta_0) + \frac{1}{2}(\hat{\beta}-\beta_0)^T
\bigl(X^TX\bigr) (\hat{\beta}-\beta_0)
\\
& \geq & - \bigl\|X^T\varepsilon\bigr\|_2 \|\hat{\beta}-
\beta_0\|_2 + \frac
{\lambda_{\min}(X^TX)}{2} \|\hat{\beta}-
\beta_0\|^2_2
\end{eqnarray*}
with $\lambda_{\min}(X^TX)$ denoting the smallest eigenvalue of
$X^TX$. We see that we can take
\[
\zeta_n = \sigma^{-1/2} n^{1/2} \frac{\|X^T\varepsilon\|_2}{\lambda
_{\min}(X^TX)}.
\]


In the high-dimensional setting, of course this fails as $\lambda
_{\min}(X^TX)=0$.
What, then, can we say about
our canonical estimator
\[
\hat{\beta}_{\lambda,\XK} = \argmin_{\beta\in\real^p} \frac
{1}{2} \|y-X\beta
\|^2_2+ \lambda h_K(\beta)?
\]
Is there even a weak law of large numbers?
Without any assumptions on $K$, the answer is clearly no: if $K=\{0\}$,
then this is the original ill-posed
problem.

Under an assumption of decomposability of $K$
recent progress has been
made in providing bounds on the estimation error in (\ref
{eqcanonical}), see
\citet{decomposable}.
The notion of decomposability in \citet{decomposable} has a precise
definition which we will not dwell on here. However,
a large set of examples of decomposable penalties are penalties of the form
%
\begin{equation}
\label{eqproduct} K = \prod_{i \in{\cal I}}K_i.
\end{equation}
That is,
convex sets that can be expressed as products of convex sets lead to
decomposable
penalties.
Every example of a decomposable norm in \citet{decomposable} has this
form except the nuclear norm. For such penalties, the generalized
gradient algorithms
described in Section~\ref{secalgorithm} decompose into many smaller
subproblems. Many efficient coordinate
descent algorithms exploit this fact, see \citet{pathwise}, \citet{glmnet}.

A similar notion of decomposability we refer to as \textit{additivity}
is explored in Taylor and Tibshirani (\citeyear{tibstayloradditive}) in which case $K$ can be
expressed as $K=A+I$ with the
sum being Minkowski addition of sets. In this case, the penalty has the form
%
\begin{equation}
\label{eqadditive} h_K(\beta) = h_A(\beta) +
h_I(\beta).
\end{equation}
One concrete example of this is the $\ell_{\infty}$ ball in $\real
^p$. For $A, I$ a partition of
$\{1,\ldots, p\}$ into \textit{active} and \textit{inactive} variables, we
can write
\[
B_{\infty} = \bigl\{(u_A,0): \|u_A
\|_{\infty} \leq1 \bigr\} + \bigl\{(0,u_I): \|u_I
\|_{\infty} \leq1 \bigr\}.
\]
Any penalty of the form (\ref{eqproduct}) can be expressed in the
form (\ref{eqadditive}) in a similar fashion.
If we are allowed to introduce a linear constraint to (\ref
{eqcanonical}), then any problem with a penalty of the form (\ref
{eqadditive}) can be expressed as a problem with a penalty of the form
(\ref{eqproduct}) subject to an additional set of linear constraints.

The weak law of large numbers presented in the literature
follow a similar path to the argument above. Of course, for precise
results, the specific
penalty as well as the data generating mechanism must be more precisely
specified.

In the interest of space,
we do not pursue such precise statements here.
Rather, we will just attempt to paraphrase these results, of which
there exist many in the literature [cf. \citet
{decomposable}, \citet{bickelritov}, Obozinski, Wainwright and
  Jordan (\citet{supportrecovery}),
\citet{highdimstat}] with \citet
{decomposable} being a
particularly nice place to read in detail. Under various assumptions on
the tails of
$\varepsilon$, as well as the assumption $h_I(\beta_0)=0$ and $K^{\circ
}$ is bounded,\footnote{Recall the definition
of the polar body of $K$
%
\begin{equation}
\label{eqpolar} K^{\circ} = \bigl\{\nu\in\real^p:
u^T\nu\leq1, \forall u \in K \bigr\}.
\end{equation}
For any $K$, the seminorm $h_{K^{\circ}}$ is the dual seminorm of $h_K$.}
the canonical result has the form for $\lambda\geq C_1 \cdot\Ee
(h_{K^{\circ}}(-X^T\varepsilon) )$
%
\begin{equation}
\label{eqerror} \|\hat{\beta}_{\lambda, \XK} - \beta_0
\|_2 \leq C_2 \cdot\sigma\frac{ \Ee(h_{K^{\circ}}(-X^T\varepsilon) ) \psi
(A)}{\kappa(A, I, X)}
\end{equation}
with high probability
for some universal $C_1, C_2$
where $\psi(A)$ is referred to as a \textit{compatibility} constant
relating the $\ell_2$ norm and the
$h_A$ seminorm; the quantity
$\kappa(A,X)$ replaces $\lambda_{\min}(X^TX)$ and is referred to as
a \textit{restricted strong
convexity} (RSC) constant.

The literature varies in their assumptions
on the noise and the design matrix $X$. For instance, in the fixed
design case one might
consider the above probability only with respect to noise, while in a
random design setting the dependence
of the constants on $X$ are typically expressed with respect to the law
that generates the design matrix $X$.

%

Having established a bound such as (\ref{eqerror}),
if one considers problems indexed by $n$, then one can obtain
a law of large numbers for the problem (\ref{eqcanonical}) so long as
the parameters
are chosen so the right-hand side decays to 0 and the bound holds with
sufficiently high probability.
Again in the interest of space, we refer
readers to the literature for more precise statements for specific
versions of the problem such as the LASSO.
We have tried to give a summary of some of the results related to the
canonical problem, though we have
barely exposed the tip of the iceberg. Under more specific assumptions
much more can be said. For example, see \citet{donohomontanari}.
%
%

\subsection{Gaussian width and metric projection: Intrinsic volumes}

For $\XK^{\circ} = X^TK^{\circ} \subset\real^n$, the error (\ref
{eqerror})
depends on the
quantity
%
\begin{equation}
\label{eqwidth} \Ee\bigl(h_{K^{\circ}}\bigl(-X^T\varepsilon\bigr)
\bigr) = \Ee\bigl(h_{\XK}(-\varepsilon)\bigr).
\end{equation}
For fixed $X$, this quantity is referred to as
the Gaussian width of $\XK$ and it also intimately
related to our first tool in the toolbox, the metric projection. Specifically,
consider the tube of radius $r$ around $\XK^{\circ}$. That is,
%
\begin{equation}
\label{eqtube} \tube\bigl(\XK^{\circ},r\bigr) = \bigl\{z \in
\real^n: \bigl\|z - \pi_{\XK
^{\circ}}(z)\bigr\|_2 \leq r \bigr\}.
\end{equation}
Then, a classical result of Steiner in the case of convex bodies and
Weyl in the case of manifolds says that
the Lebesgue measure of the tube, assuming $\XK^{\circ}$ is bounded,
can be expressed as
%
\begin{equation}
\bigl\llvert\tube\bigl(\XK^{\circ},r\bigr)\bigr\rrvert_{\real^n} =
\sum_{j=0}^n r^j
\omega_j \lips_{n-j}\bigl(\XK^{\circ}\bigr), \qquad r
\leq r_c\bigl(\XK^{\circ
}\bigr)=\infty,
\end{equation}
where the $\lips_l(\XK^{\circ})$ are referred to as the intrinsic
volumes of $\XK^{\circ}$
and $\omega_l = |B_2(1)|_{\real^l}$ is the Lebesgue measure of the
unit $\ell_2$ ball in $\real^l$. See \citet
{rfg}, \citet{federer}, \citet{schneider}, \citet{weyl} for more details on such volume of tubes formulae.
When $\XK^{\circ}$ is unbounded, Federer's curvature measures
[\citet
{federer}] can be used to define the volume of local tubular neighborhoods.
Using Gaussian process techniques [\citet{vitale}], it can be
shown that
%
\begin{equation}
\label{eqgausswidth} \lips_1\bigl(\XK^{\circ}\bigr) = (2
\pi)^{-1/2}\Ee\bigl(h_{\XK^{\circ}}(\varepsilon)|X\bigr), 
\end{equation}
where $\varepsilon|X \sim N(0, I_{n \times n})$.

For $D=\XK$ a smooth domain, that is convex set with non-empty
interior bounded by a smooth hypersurface
and for $j \leq n-1$
%
\begin{equation}
\label{eqintrinsicvolumes} \lips_j(D) \propto\int
_{\partial D}P_{p-j-1}(\lambda_{1,x},\ldots,
\lambda_{p-1,x}) \operatorname{Vol}_{\partial D}(dx)
\end{equation}
with $P_j$ the $j$th elementary symmetric polynomial of the so-called
principal curvatures of $\partial D$ at $x$ [\citet{rfg}, \citet{noniso}]. These
are just the eigenvalues of the second fundamental form in the unit
inward normal direction. This formula can be derived
by considering the inverse of the metric projection map (\ref
{eqmetricproj}). The inverse takes $(x, \eta_x)$ defined on the
extended outward normal bundle of
$K^{\circ}$. The inverse of the map is simply the exponential map
restricted to the outward normal bundle, or, more simply
%
\begin{equation}
\label{eqexponential} x \mapsto x + \eta_x.
\end{equation}
A fairly straightforward calculation
yields the relation (\ref{eqintrinsicvolumes}).

The main take away message above is that the functionals in the tube
formula, such as the Gaussian width
(\ref{eqwidth}), are related to properties of the metric projection
map onto~$\XK$.

As written above, the intrinsic volumes are defined implicitly through
a volume calculation, and it is not clear that they
extend to the infinite dimensional setting. Under the right conditions
of course, such extensions are indeed possible. See \citet{vitale} for
a nice discussion of this problem. An alternative definition of
intrinsic volumes specific to the Gaussian case
was considered in \citet{taylorvadlamani}, which were defined as
coefficients in an expansion of the \textit{Gaussian} measure of
$\tube(\XK^{\circ},r)$ as described in the Gaussian Kinematic
Formula [\citet{rfg}, \citet{taylorkinematic}].

\subsection{Risk estimation}

Another quantity of interest for our problem (\ref{eqcanonical}) is
an estimate of
how much ``fitting'' we are performing, as a function of $\lambda$.
One quantitative
measure of this is captured by Stein's estimate of risk [\citet
{stein}], also known as
SURE. Suppose now that $y \sim N(\mu, I)$ and we estimate $\mu$ by
the estimator
$\hat{\mu}$. The SURE estimate is an unbiased estimate of
\[
\operatorname{Risk}(\hat{\mu}) = \Ee\bigl(\|\mu- \hat{\mu}\|^2_2
\bigr).
\]
The estimated degrees of freedom of this estimator is one part of the
SURE estimate and is defined as
%
\begin{equation}
\widehat{\operatorname{Cov}(y, \hat{\mu}}) = \operatorname{div}\bigl
(\nabla
\hat{\mu}(y)\bigr).
\end{equation}
%

Suppose $y \sim N(\mu, I)$ and consider our residual form of the
original estimation problem for $\mu=\Ee(y)$
%
\begin{equation}
\hat{\mu}(y) = y - \pi_{\lambda\XK}(y).
\end{equation}
If $K$ possesses a nice stratification, as do all the examples
mentioned above, then, for almost every $y$, $\pi_{\lambda\XK}(y)$
is in the relative interior\vadjust{\goodbreak}
of some fixed stratum ${\cal S}$ of the normal bundle of $\XK$ over
which the dimension
of the tangent space is constant, and the normal bundle has a locally
conic structure
${\cal S} = {\cal T} \times{\cal N}$ of tangent and normal directions
[\citet{schneider}, \citet{rfg}].
Having fixed this stratum, we can write
$y=x+\eta_x$ in (ortho)normal coordinates centered at $(y-\hat{\mu
}(y), \hat{\mu}(y))$. In these
coordinates
%
\begin{equation}
\hat{\mu}\bigl(y(x,\eta_x)\bigr) = \eta_x.
\end{equation}
In order to relate the above to the problem (\ref{eqcanonical}), one should
invert the above chart to find $(x,\eta_x)$ in terms of $y$.
In the residual form (\ref{eqresid}) it is easy to show that
\begin{eqnarray*}
\bigl(x(y), \eta_x(y) \bigr) &=& \biggl(
\argmin_{r \in\lambda
\cdot\XK} \frac{1}{2} \|y-r\|^2_2, y -
\argmin_{r \in\lambda\cdot
\XK} \frac{1}{2} \|y-r\|^2_2
\biggr)
\\
&=& \bigl(y - \hat{\mu}(y), \hat{\mu}(y)\bigr).
\end{eqnarray*}
The derivatives along the normal directions yield
a purely dimensional term, while directions in the tangent directions yield
curvature terms. This observation is enough to derive the following
form of the degrees of freedom
%
\begin{equation}
\label{eqSURE}
\operatorname{div}\bigl(\nabla\hat{
\mu}(y)\bigr)
= n - \operatorname{dim}({\cal
T}_y) + \operatorname{Tr}(S_{(y - \hat
{\mu}(y),\hat{\mu}(y))}).
\end{equation}
Above, ${\cal T}_y$ is the tangential part of the stratum containing
$x(y)$ and $S_{(x,\mu_x)}$ is the second fundamental form of
${\cal T}_y$ in $\lambda\cdot\XK$ as described [\citet{rfg}]. When
$K$ is a polyhedral set, the second term disappears
and the degrees of freedom can be computed by computing the rank of a
certain matrix
[\citet{tibshiranidegrees2012}, \citet{bienhierarchical}].
%

\subsection{Hypothesis testing: Weak convergence for high-dimensional
inference?}

Another fundamental tool in inference for least squares models is the
ability to form hypothesis
tests, as well as confidence intervals for the ``true'' mean. Such concepts
clearly need a model, which we might take to be the usual model
\[
y \sim N(\mu,I).
\]
Under the assumption that $\mu=X\beta_0$, classical inference in
linear models (assuming $X^TX$ is full rank) yields confidence
intervals and hypothesis tests for any linear functional $\nu^T\mu$
based on the
coefficients $\beta_0$.

What can we say about our canonical problem (\ref{eqcanonical})? This
is an area in which we still don't know all the answers. In some sense,
we are in the situation analogous to Bernoulli having proved a weak
law of large numbers without the central limit theorem.
Some progress has been made for specific models of the design matrix
$X$ for the LASSO as well as group LASSO models, see
\citet{meinshausenbuhlmann}, \citet{buhlmann}, \citet{minniercai},
\citet{murphylaber}, Wasserman and Roeder (\citeyear{wassermanroeder}), \citet{donohomontanari}.\vadjust{\goodbreak}

Other recent work [\citet{significancelasso}] gives some hints at what
inferential tools may prove useful in this weak convergence theory. As
described in the LARS algorithm
[\citet{lars}, \citet{tibshiraniunique}], an entire path of solutions $\hat
{\beta}_{\lambda,\XK}$ can be formed for the LASSO, that is, when
$K=B_{\infty}$. These paths are piecewise linear, with knots at points
where the \textit{active set} changes. The covariance statistic
measures some change in the correlation of the fitted values between
two knots $\lambda_k$ and $\lambda_{k+1}$ in the LASSO path. It has
the form
%
\begin{equation}
\label{eqcovstat} T_k = C_k \lambda_k (
\lambda_k - \lambda_{k+1})
\end{equation}
for some random scaling $C_k$ related to the active set and the
variable added to the active set at $\lambda_{k+1}$.
The form for $k=1$ is particularly simple: suppose that $\hat{j}_1$
and $X$ is such the first variable in the LARS path, that is,
\[
\hat{j}_1 \in\argmax_j \bigl|X_j^Ty\bigr|.
\]
Then,
%
\begin{equation}
\label{eqcovstat1} T_1 = \frac{1}{\sigma^2} y^TX\hat{
\beta}_{\lambda_2} = y^T X_{\hat
{j}_1} \hat{
\beta}_{\lambda_2,\hat{j}_1}.
\end{equation}
For $k \geq1$, the form of the test statistic is slightly more
complicated, though it can be expressed in terms of
$A=A_k$, the active set at step $k$ as well as $s_{A_k}$, the signs of
the active variables at step $k$ as well as
the active set $A_{k+1}$ and $s_{A_{k+1}}$, see Section 2.3 of
\citet
{significancelasso} for the full expression.
For a wide range of (sequences) covariance matrices, if the active set
at $\lambda_k$ already contains all the strong active variables, then
it is shown in \citet{significancelasso} that
%
\begin{equation}
\label{eqexp1} T_k \stackrel{D} {\rightarrow} \operatorname{Exp}(1).
\end{equation}
In particular, under the global null $y \sim N(0,\sigma^2 I)$ as long
as the design matrix satisfies some minimum growth condition,
$T_1 \stackrel{D}{\rightarrow} \operatorname{Exp}(1)$.
The main tools used in the above proof relate to the maxima of
(discrete) Gaussian processes and the generalization of an
argument previously applied to smooth Gaussian processes in \citet
{taylorvalidity} and \citet{rfg}.
Ongoing work suggests that such a limiting distribution will hold for
many (sequences) of $K$ and design matrices $X$.

As for confidence intervals for the parameters related to $k$ strong
variables, the relation to extreme values of Gaussian processes suggest
that the bias-corrected or relaxed LASSO estimate of the active
coefficients will have accurate coverage. This is ongoing work.

\section{Smooth Gaussian processes: Relaxing convexity}
\label{secgaussian}

In the special case that $K$ is a cone, we saw that the distribution of
a particular
likelihood ratio test could be expressed in terms of the supremum of a
Gaussian process indexed by a subset of the sphere.
Equivalently, this supremum could be expressed in terms of the metric
projection onto the cone
$\XK$. It is well-known that the distribution\vadjust{\goodbreak} of this likelihood ratio
test statistic is a mixture
of $\chi^2$'s of varying degrees of freedom. This distribution is
sometimes referred to as a $\bar{\chi}_2$ distribution.

The mixture weights can be expressed in terms of the geometry of the
set $M = S(X,K) = \XK\cap S(\real^n)$. In particular, it is known
[\citet{rfg}, Takemura and Kuriki
(\citeyear{takemurakuriki1,takemurakuriki2})] that if $\varepsilon\sim
N(0, I)$ for $u > 0$
%
\begingroup
\abovedisplayskip=6.5pt
\belowdisplayskip=6.5pt
\begin{equation}
\label{eqsup} \Pp\Bigl(\sup_{\nu\in M}\bigl( \nu^T
\varepsilon\bigr)^+ > u \Bigr) = \sum_{j=0}^{\infty}
\lips_j(M) \rho_j(u),
\end{equation}
where
%
\begin{equation}
\label{eqrho} \rho_j(u) = \cases{\displaystyle \int_u^{\infty}
\frac{1}{\sqrt{2\pi}} e^{-t^2/2} \,dt, & $j = 0$,
\vspace*{2pt}\cr
\displaystyle (2\pi)^{-(j+1)/2}
H_{j-1}(u) e^{-u^2/2}, & $j \geq1$,}
\end{equation}
and
\[
H_j(u) =(-1)^j \frac{\partial^j}{\partial u^j} e^{-u^2/2}
\]
are the standard Hermite polynomials. The functions (\ref{eqrho})
are known as the EC or Euler characteristic densities for a
Gaussian field, see \citet
{rfg}, \citet{takemurakuriki2}, \citet{worsleyboundary},
\citet{worsleyunified1}, \citet{worsleyunified2}.
While the sum above is written as an infinite sum it terminates at
$\operatorname{dim}(M)$.

Note that this implies
\[
\Ee\Bigl(\sup_{\nu\in M} \bigl(\nu^T \varepsilon\bigr)^+
\Bigr) = \frac
{1}{\sqrt{8\pi}} \lips_1(M)
\]
which is an analogous way to derive Gaussian width (\ref{eqgausswidth}).

One of the derivations of the above formula, the so-called \textit{volume
of tubes} approach
[\citet{siegmund}, \citet{sun}, \citet{takemurakuriki1}] involves studying the Jacobian
of the inverse of the
spherical metric projection map (\ref{eqmetricprojsphere}), that
is, the exponential
map on $S(\real^n)$ which sends a pair $(x,\eta_x)$ to
\[
\cos\bigl(\|\eta_x\|\bigr) \cdot x + \sin\bigl(\|\eta_x\|\bigr) \cdot
\eta_x.
\]
Another approach, the expected Euler characteristic approach
[\citet{rfg}, \citet{worsleyboundary}] involves counting critical points of
the Gaussian process above the level $u$ according to saddle type and applying
Morse theory and the Rice--Kac formula [\citet{azais}, \citet{rfg}] to count the
expected number of such points.

Neither of these approaches strictly require convexity of the cone
generated by the parameter set $M$. Rather,
they depend on a notion of local or infinitesimal convexity referred to
as positive
reach [\citet{federer}]. Hence, the parameter sets may have finite
critical radius.
They are both approaches used to form an approximation
%
\begin{equation}
\label{eqsupgeneric} \Pp\Bigl(\sup_{\nu\in M} f(\nu) > u \Bigr)
\approxeq\sum_{j=0}^{\infty}
\lips_j(M) \rho_j(u)\vadjust{\goodbreak}
\end{equation}
\endgroup
for some centered, smooth, Gaussian process $f$ having constant
variance 1 on a
(possibly stratified) manifold $M$.
In the volume of tubes approach, the critical radius appears in a
natural way and enters into an estimate of the error of the volume of
tubes approach. In both approaches, though it is clearer in the
expected Euler characteristic approach, the spherical critical radius
is in fact the spherical critical radius of
$\Psi(M)$ where $\Psi:M \rightarrow S(H)$ where $S(H)$ is the unit
sphere in $H$, the reproducing kernel Hilbert space of
$f$. Either approach yields roughly the same estimate of error: for $u$
large enough
%
\begin{equation}
\label{eqsupgenerictail} \liminf_{u \rightarrow\infty} -\frac{2}{u^2}
\Biggl\llvert\Pp\Bigl(\sup_{\nu\in M} f(\nu) > u \Bigr) - \sum
_{j=0}^{\infty} \lips_j(M)
\rho_j(u) \Biggr\rrvert= 1 + \tan^2 \bigl(r_c(M)
\bigr).
\end{equation}
The above says that the \textit{relative error} in the approximation is
exponentially small whenever $r_c(M) > 0$.

The critical radius of $M$ in the expected Euler characteristic
approach arises in terms of
a functional of the original process $f$. Specifically, if we assume
$M$ is a manifold without boundary, then define for each $x\neq y$ the
process introduced in
\citet{taylorvalidity}
%
\begin{equation}
\label{eqfxy} f^x(y)= \frac{f(y) - \Ee(f(y)|f(x), \nabla f(x))}{1 - \Ee
(f(x) \cdot f(y))}.
\end{equation}
Then,
%
\begin{equation}
\label{eqfxy2} \cot^2\bigl(r_c(M)\bigr) = \sup
_{x \neq y} \Ee\bigl(f^x(y)^2\bigr).
\end{equation}
Hence, the critical radius depends in an explicit way on the covariance
function of the
process $f$.

As mentioned previously, the argument related to derivation of (\ref
{eqfxy2}) in the smooth case led directly to the
exponential limit in (\ref{eqexp1}). Such a connection suggests a
relation between the
distribution of the maxima of smooth Gaussian processes, specifically
the spacings of the extreme values, can be
used to derive weak convergence results for high-dimensional inference.
This is ongoing work.
%
%
%

\section{Conclusion}

We have described what we call the two basic tools of the geometry of
least squares that are just as relevant as when Gauss and Legendre
disputed their original discovery over 200 years ago. While these tools
are not the most technically sophisticated tools, ceding that ground to
exponential families for the canonical model (\ref{eqcanonical}) and
empirical processes for the fluctuation theory in
Section~\ref{secinference}, they nevertheless provide guiding
principles for these more precise tools. We would argue that the
Gaussian picture provided by the geometry of least squares, gets much
of the picture correct under sufficient moment conditions. For heavier
tailed results, of course much of Section~\ref{secinference} would have
to be reframed and Section~\ref{secgaussian} paints quite a different
picture [Adler, Samorodnitsky and Taylor
(\citeyear{samorodnitsky1,samorodnitsky2})].

As Bernoulli found himself with just a law of large numbers, the field
of statistics is roughly at this same stage in high-dimensional
inference. We are hopeful that the geometry of least squares will
eventually guide the field to weak convergence results in
high-dimensional inference for the canonical problem
(\ref{eqcanonical}).

\section*{Acknowledgement}

Supported in part by NSF Grant DMS-12-08857 and
AFOSR Grant 113039.


%

\printhistory

\end{document}